\documentclass[12pt,a4paper,oneside,onecolumn]{article}
\title{On Gel'fand-Kolmogoroff type results}
\author{Zihindula Mushengezi E.}
\usepackage{etex}
\usepackage[]{latexsym}
\usepackage[T1]{fontenc}
\usepackage[latin1]{inputenc}
\usepackage[english]{babel}
\usepackage{graphics}
\usepackage[]{amssymb}
\usepackage{multicol}
\usepackage{amsmath}
\usepackage{graphicx}
\usepackage{amsfonts}
\usepackage[all]{xy}
\usepackage{pifont}
\usepackage{microtype}
\usepackage{answers}
\usepackage{tabularx}
\usepackage{float}
\usepackage{color}

\usepackage{listings}
\definecolor{darkWhite}{rgb}{0.94,0.94,0.94}
\usepackage{palatino}

\newcount\exenum     \newdimen\numindent
\newcount\qqnum     \newdimen\qqmarge
\newcount\sqnum     \newdimen\sqmarge

\outer\def\bye{%
\vskip 0pt\@endmulticol\@endgroup              
\ifanswer \let\next=\exobye@                   
\else     \let\next=\@exobye                   
\fi\next}

\usepackage[Glenn]{fncychap}

\usepackage{pstricks}
\usepackage{pgf,tikz}
\usepackage{pstricks-add}
\usepackage{pst-plot}
\usetikzlibrary{arrows}
\usepackage{pgfplots}
\usepackage{pgfkeys}

\definecolor{qqqqqq}{rgb}{0,0,0}

\definecolor{xdxdff}{rgb}{0.49,0.49,1}
\definecolor{qqwuqq}{rgb}{0,0.39,0}
\definecolor{qqqqff}{rgb}{0,0,1}
\definecolor{ttttff}{rgb}{0.2,0.2,1}
\definecolor{uququq}{rgb}{0.25,0.25,0.25}

\renewcommand{\epsilon}{ \varepsilon}

\newcommand{\euro}{\texteuro{}}

\newcommand{\pre}{{\bf Proof.\ }}

\newcommand{\ve}{{\bf e}} 
\newcommand{\vu}{{\bf u}}

\newcommand{\vv}{{\bf v}}

\newcommand{\vy}{{\bf y}}

\newcommand{\vh}{{\bf h}}
\newcommand{\vk}{{\bf k}}

\newcommand{\vU}{{\bf U}}

\newcommand{\cl }{\mathcal }

\newtheorem{theo}{Theorem}[section]
\newtheorem{prop}[theo]{Proposition}
\newtheorem{lem}[theo]{Lemma}
\newtheorem{cor}[theo]{Corollary}
\newtheorem{rem}[theo]{Remark}

\newcounter{exercice}

\definecolor{fbase}{rgb}{0.8,0.8,1}
\definecolor{fgris}{gray}{0.6}
\definecolor{frouge}{HTML}{DC143C}
\definecolor{fvert}{rgb}{0.6,1,0.6}
\definecolor{fbleu}{rgb}{0.4,0.4,1}
\definecolor{fjaune}{HTML}{DCDC14}
{\endMakeFramed}
\makeatletter
\DeclareRobustCommand\sfrac[1]{\@ifnextchar/{\@sfrac{#1}}%
                                            {\@sfrac{#1}/}}
\def\@sfrac#1/#2{\leavevmode\kern.1em\raise.5ex
         \hbox{$\m@th{\fontsize\sf@size\z@\selectfont#1}$}
         \kern-.1em/\kern-.15em\lower.55ex
          \hbox{$\m@th{\fontsize\sf@size\z@\selectfont#2}$}}

\DeclareRobustCommand{\Efrac}[2]{{\displaystyle\begingroup
\raise2ex\hbox{$\m@th{#1}$}\endgroup\@@over \lower1ex
\hbox{$\m@th{#2}$}}}
\makeatother

\newtheorem{theorem}{Theorem}[section]

\numberwithin{equation}{section}

\newtheorem{Exc}{Exercice}
\Newassociation{correction}{Soln}{mycor}
\Newassociation{indication}{Indi}{myind}

\def\exo#1{\futurelet\testchar\MaybeOptArgmyexoo}
\def\MaybeOptArgmyexoo{\ifx[\testchar \let\next\OptArgmyexoo
                        \else \let\next\NoOptArgmyexoo \fi \next}
\def\OptArgmyexoo[#1]{\begin{Exc}[#1]\normalfont}
\def\NoOptArgmyexoo{\begin{Exc}\normalfont}

\newcommand{\finexo}{\end{Exc}}
\newcommand{\flag}[1]{}

\tikzset{
xmin/.store in=\xmin, xmin/.default=-3, xmin=-3,
xmax/.store in=\xmax, xmax/.default=3, xmax=3,
ymin/.store in=\ymin, ymin/.default=-3, ymin=-3,
ymax/.store in=\ymax, ymax/.default=3, ymax=3,
}

\newcommand{\entete}[1]

\begin{document}
 \Opensolutionfile{mycor}[ficcorex]
 \Opensolutionfile{myind}[ficind]
 \entete{\'Enoncés}


\maketitle

\begin{abstract}
We prove that a vector bundle $ E \to M$ is characterized by the associative structure of the space of symbols of the Lie algebra generated by all differential operators on $E$ which are eigenvectors of the Lie derivative in the direction of the Euler vector field.

We also obtain similar result with the $\mathbb{R}-$ algebra of smooth functions which are polynomial along the fibers of $E.$
This allows us to deduce a  Gel'fand-Kolmogoroff type result for the $\mathbb{R}-$algebra  ${\rm Pol}(T^*(M))$ of symbols of the differential operators of $M.$
\end{abstract}

\section*{Introduction}

An associative algebra $\cl A(M)$ characterizes a topological space $M$ if, for any topological space $N$; the associative algebras $\cl A(M)$ and $\cl A(N)$ are isomorphic if, and only if, the topological spaces $M$ and $N$ are homeomorphic.

The algebraic characterization of topological spaces dates back to the end of the 30s  of the last century with I. Gel'fand and A. Kolmogoroff in : \cite{GelKol}.\\ 
It is established in this paper that a compact topological space $ M $ is characterized by the associative algebra $ {\rm C} (M) $ of continuous functions on $M$ with  values   in   $\mathbb{R} $ or $ \mathbb{C}. $

The methods developed in \cite{GelKol} can be applied to the case of a smooth manifold $ M, $ assumed of Hausdorff and with a countable basis, by considering the algebra $ \cl A(M) = {\rm C}^ \infty (M) $ of smooth functions on $ M. $ 

Any isomorphism of  $\mathbb{R}-$algebras $ \Psi: {\rm C}^ \infty (M) \to {\rm C}^ \infty (N) $ has the form
 \[
   \Psi: f \mapsto f \circ \psi,
 \]
with $ \psi: N \to M $ a diffeomorphism of smooth manifolds. This result is due to Milnor.

The Lie-algebraic characterization of differential manifolds appeared for the first time in 1954 in a paper of P.E. Pursell and M.E. Shanks. More precisely in \cite{PurShan}.\\

One way to generalize these results is to establish characterizations with Lie algebras larger than that of vector fields of a manifold. This is what J. Grabowski and N. Poncin proposed in \cite{GraPon2, GraPon1, GraPon4}.\\
We can also obtain a generalization of the result of Pursell and Shanks by characterizing vector bundles by Lie algebras. In this way, let us cite \cite{Lec1,LecLeuZih,LecZih1,LecZih2}.\\ 
In \cite{Lec1} it is established that, under certain assumptions, a vector bundle can be characterized by the Lie algebra of its infinitesimal automorphisms.\\

In this work, we generalize the initial result of Gel'fand and Kolmogoroff in two ways. On the one hand, by considering much larger associative algebras, and on the other, by characterizing  vector bundles rather than manifolds.\\
This generalization was made possible thanks to the algebraic aspects presented in \cite{Jet} by Nestruev J.

\section{From results due to Grabowski and Poncin}\label{GraPo}

Let $ M $ be a manifold. We know that, thanks to the work of Grabowski J. and Poncin N., the algebras $ \cl D (M),$ of all differential operators of $M$, and $ \cl S (M),$ of the symbols of these operators, characterize the manifold $ M $ by their respective structures of Lie algebras. 
For their structures of $ \mathbb{R}-$ associative algebras, the following result provides a partial answer to the question.
\begin{theo}\label{caract alg variétés}
Let $M$ and $N$ be smooth manifolds.\\
The $ \mathbb{R}-$ associative algebras $ \cl D (M) $ and $ \cl D (N) $ are isomorphic if, and only if, the smooth manifolds $ M $ and $ N $ are diffeomorphic.
\end{theo}
\pre
The assertion comes from the fact that any isomorphism of associative algebras between these algebras preserves their Lie bracket, the latter being nothing other than the bracket of commutators. \hfill$\blacksquare$\\

For the Poisson algebra $\cl S(M),$ obtaining the characterization of the manifold $ M $ by its associative algebra structure requires another theoretical approach.\\
This is precisely what we propose in the next section.

\section{Some notions of algebra}

In this part, we present some notions on the $ \mathbb {R} - $commutative associative algebras which we use in the following, to obtain results of the Gel'fand-Kolmogoroff type. They are mainly taken from \cite{Jet}.\\
Let $ \cl F $ be a $ \mathbb{R}-$ commutative associative algebra with unit. In the following lines, we will simply say that $ \cl F $ is a $ \mathbb{R}-$ algebra. 
 We denote by $ | \cl F |: = M, $ the set of all homomorphisms of $ \mathbb{R}-$ algebras of $\cl F $ in $\mathbb{R}$
\[
  M\ni x:\cl F\to \mathbb{R}: f\mapsto x(f).
\]
The elements of $ M $ are then called \textit{$ \mathbb{R} - $ points} of the algebra $ \cl F $ and the set $ | \cl F |, $ the \textit{ dual space } of $ \mathbb{R} - $ points. \\
A $ \mathbb {R} - $ algebra $ \cl F $ is said to be \textit{geometric} if the subset
\[
  \cl I(\cl F)=\bigcap_{x\in|\cl F|}Ker\, x
\]
contains only the null element.
\begin{prop}\label{alg des fx sur ens}
If $\cl F$ is a $\mathbb{R}-$algebra of functions over a given set $ N, $ then $ \cl F $ is geometric.
\end{prop}
\pre
Let us specify that by function on a set $ N $ we mean a function $ f: N \to \mathbb{R}. $
Consider the application $ \theta: N \to | \cl F | $ associating to any point $ a \in N, $ the homomorphism of $\mathbb{R}-$algebras $\cl F\to \mathbb{R}:f\mapsto f(a).$ In other words, 
\[
  \theta(a)(f)=f(a) ,\quad\forall a\in N,\forall f\in\cl F,
\]
and we have that $ \theta (a) $ is a homomorphism of $ \mathbb {R} - $ algebras.
Observe that through the following definition
\[
f(x):=x(f), x\in|\cl F|, \, f\in\cl F,
\]
any element of $ \cl F $ can be seen as a function on the dual space $ | \cl F |. $ In particular, we have
\[
 f(\theta(a))=\theta(a)(f)=f(a).
\]
Therefore, if an element $ f \in \cl F $ vanishes in $ \theta (a), $ for all $ a \in N, $ then $ f $ is the null element of the $ \mathbb { R} - $ algebra $ \cl F.$
From the inclusion
\[
  (\bigcap_{x\in|\cl F|}Ker x)\subset (\bigcap_{a\in N}Ker \theta(a))
\]
we then deduce that $ \cl F $ is geometric. \hfill$\blacksquare$

A $ \mathbb {R} - $ geometric algebra $ \cl F $ is called \textit{$ {\rm C}^\infty- $ closed} if for any finite collection of elements $ f_1, \cdots, f_k \in \cl F $ and any function $ g \in {\rm C}^\infty (\mathbb {R}^k), $ there exists $ f \in \cl F $ such that 
 \[
  f(a)=g(f_1(a),\cdots,f_k(a)),\quad \forall a\in|\cl F|.
 \]
Note that such an element $ f \in \cl F $ is uniquely determined by $ f_i $ and $ g $ because $ \cl F $ is geometric. 

As an example, for any $ n \in \mathbb {N} _0, $ the $ \mathbb {R} - $ algebra ${\rm C}^\infty(\mathbb{R}^n)$ is ${\rm C}^\infty-$closed.\\
 
Let $ \cl F $ be a $ \mathbb {R} - $ geometric algebra that we identify with a subset of the algebra of functions over $ | \cl F |. $ Consider the set $ \overline {\cl F} $ defined as the set of functions over $ | \cl F | $ that can be written in the form
 \[
   g(f_1,\cdots,f_k),\quad k\in \mathbb{N},\,\,f_i\in \cl F,\,\, 
                                                g\in{\rm C}^\infty(\mathbb{R}^k).
 \]
The set $ \overline {\cl F} $ has an obvious structure of $ \mathbb {R} - $ algebra and $ \cl F $ is indeed a subalgebra.
By virtue of the previous proposition \ref{alg des fx sur ens}, $ \overline {\cl F} $ is geometric.
This algebra is also $ {\rm C} ^\infty- $ closed.\\

By definition, the $ \mathbb {R} - $ algebra $ \overline {\cl F} $ thus associated with the $ \mathbb {R} - $ geometric algebra $ \cl F $ is called \textit{smooth envelope} of $ \cl F. $
It has the following significant property.
\begin{prop}\label{prop univer c infi fermé}
Let $ \cl F $ be a $ \mathbb {R} - $ geometric algebra and $ \overline {\cl F} $ its smooth envelope. For any homomorphism $ \Psi: \cl F \to \cl F'$ of\,\, $ \mathbb{R} - $ algebras, with $ \cl F' $ a $ \mathbb {R} - $ algebra $ {\rm C }^\infty- $ closed, there exists a unique homomorphism of\,\, $ \mathbb {R} - $ algebras $ \overline {\Psi}: \overline {\cl F} \to \cl F '$ such that $ \Psi = \overline {\Psi} \circ i, $ with $ i: \cl F \to \overline {\cl F} $ the canonical inclusion.
\end{prop}
\pre
The following map
\[
  \overline{\Psi}(g(f_1,\cdots,f_k))=g(\Psi(f_1),\cdots,\Psi(f_k))
\]
defined from $ \overline {\cl F} $ to $ \cl F ',$ satisfies the conditions of the proposition.\hfill$\blacksquare$\\

We then deduce the following result for which a proof can be found in \cite{Leu}.
\begin{prop}\label{Thomas env lisse}
 Let $ \cl F_1 $ and $ \cl F_2 $ be two\, $ \mathbb {R} - $ geometric algebras. If $ \Psi: \cl F_1 \to \cl F_2 $ is an isomorphism of\,\, $ \mathbb {R} - $ algebras, there then exists a unique isomorphism of\,\, $\mathbb {R} - $ algebras $ \widehat{\Psi }: \overline {\cl F_1} \to \overline {\cl F_2} $ such that $ \Psi =\widehat{\Psi} \circ i_1, $ with $ i_1: \cl F_1 \to \overline {\cl F_1} $ the canonical inclusion.
\end{prop}
\begin{prop}\label{inclusion egalite envpe lisse}
Let $ \cl F $ and $ \cl G $ be two $ \mathbb {R} - $ geometric algebras. We then have the following implication
 \[
   \cl F\subset \overline{\cl G}\Longrightarrow \overline{\cl F}\subset \overline{\cl G}.
 \]
And one deduces the below equality
  \[
    \overline{\cl F}=\overline{\overline{\cl F}}\cdot
  \]
\end{prop}
\pre
As $ \cl F \subset \overline {\cl G}, $ denote by $ i: \cl F \to \overline {\cl G} $ the canonical injection. According to Proposition \ref{prop univer c infi fermé}, $ i $ extends into a homomorphism of $ \mathbb {R} - $ algebras $ \overline {i}: \overline {\cl F} \to \overline{\cl G} $. The formula defining $ \overline {i} $ given at the beginning of the proof of the same Proposition \ref{prop univer c infi fermé} shows that $ \overline {i} $ is nothing other than the canonical injection of $\overline{\cl F}$ in $\overline{\cl G}.$ \\
From the inclusion $\overline{\cl F}\subset\overline{\cl F},$ we then deduce that $\overline{\overline{\cl F}}\subset \overline{\cl F}.$ Hence the conclusion.\hfill $\blacksquare$

Let us state a last result in this section and whose proof is given in \cite{Jet}.
\begin{prop}\label{prop Jet env lisse des polynomes}
Let $ E \to M $ be a vector bundle. The smooth envelope of the $ \mathbb {R} - $ algebra $ \cl A (E) $ of fiberwise polynomial functions of $E$ is the $\mathbb{R}-$algebra ${\rm C}^\infty(E).$
\end{prop}

\section{A Gel'fand-Kolmogoroff type result}
We establish in the following lines an algebraic characterization of a vector bundle with its $\mathbb{R}-$ algebra of fiberwise polynomial functions. 
\begin{lem}\label{lem:degree 0}
Let $\pi : E \to M$ and $\eta : F \to N$ be two vector bundles.\\ 
If $\Psi: {{\mathcal A}}(E) \to {{\mathcal A}}(F)$ is an isomorphism of associative algebras, then
 \[
  \Psi({{\mathcal A}}^0(E)) = {{\mathcal A}}^0(F)
 \]
\end{lem}
\pre
For any nonvanishing element $u \in {\mathcal A}^0(E)$, the function $u^{-1} : e \mapsto
\frac{1}{u(e)}$ is still an element of ${\mathcal A}^0(E)$. Since $\Psi$ is an homomorphism, we have
\[
{\mathcal A}^0(F) \ni 1_F = \Psi(1_E) = \Psi(u . u^{-1}) = \Psi(u) . \Psi(u^{-1}) 
\]
This implies that $\Psi(u)$ and $\Psi(u^{-1})$ are two nonvanishing functions, fiberwise polynomial with zero degree. 
Since the function  $u^2 + 1_{E}$ is nonvanishing and belongs to ${\mathcal A}^0(E),$ for any element $u \in {\mathcal A}^0(E)$, it follows that 
\[
{\mathcal A}^0(F) \ni \Psi(u^2+1_{E}) = \Psi(u).\Psi(u) + 1_{F}\cdot
\]
This shows that $\Psi(u)$ is fiberwise polynomial and of zero degree.
We just proved inclusion $\Psi({{\mathcal A}}^0(E)) \subset {{\mathcal A}}^0(F)$. The conclusion
follows by applying the same reasoning to the inverse homomorphism $\Psi^{-1}$. \hfill $\blacksquare$
\begin{lem}\label{lem:iso poly induit iso linéaire}
  Any  diffeomorphism $\varphi:\mathbb{R}^n \to \mathbb{R}^n$  such that
    \[
      {\rm Pol }(\mathbb{R}^n)=\{P\circ\varphi: P\in \rm{ Pol }(\mathbb{R}^n)\}
    \]
  is polynomial and its " linear part" coincides with the automorphism  $\varphi_{*_0}$ of the $\mathbb{R}-$vector space  $\mathbb{R}^n.$  
\end{lem}
\pre
Let us first observe that $ \varphi $ is necessarily polynomial; in the sense that we can write
 \[
   \varphi(y_1,\ldots,y_n)=(\varphi_1(y_1,\ldots,y_n),\varphi_2(y_1,\ldots,y_n),\cdots,\varphi_n(y_1,\ldots,y_n))
 \]
 with, for any  $j\in\{1,\ldots,n\},$
  \[
   \varphi_j(y_1,\ldots,y_n)=\lambda_0^j+\lambda_1^jy_1+\cdots+\lambda_n^jy_n+\varphi_{\geq 2}^j(y_1,\ldots,y_n)\cdot
  \]
Note that the expression $\varphi_{\geq 2}^j(y_1,\ldots,y_n)$ is polynomial in $y_1,\ldots,y_n$ 
of degree greater than or equal to 2. Indeed, this comes from the fact that, for any $j\in\{1,\ldots,n\},$ the map
  \[
    (y_1,\ldots,y_n)\in\mathbb{R}^n\mapsto y_j\in\mathbb{R}
  \]
is an element of $\rm{ Pol }(\mathbb{R}^n).$
We have that for any $\vy=(y_1,\ldots,y_n)\in\mathbb{R}^n,$
 \[
  \varphi_{*_\vy}= \left(
     \begin{array}{ccc}
      \lambda_1^1+\partial_1\varphi_{\geq2}^1(\vy) & \cdots & \lambda_n^1+\partial_n\varphi_{\geq2}^1(\vy)\\
      \vdots                                       &        & \vdots\\
      \lambda_1^n+\partial_1\varphi_{\geq2}^n(\vy) & \cdots & \lambda_n^n+\partial_n\varphi_{\geq2}^n(\vy)
     
     \end{array}
     \right)
 \]
is a linear automorphism of $ \mathbb {R}^n. $ In particular, for $ \vy = 0, $ we have that the matrix 
 \[
   A_\varphi=\left(
     \begin{array}{ccc}
      \lambda_1^1       & \cdots & \lambda_n^1\\
      \vdots            &        & \vdots\\
      \lambda_1^n       & \cdots & \lambda_n^n    
     \end{array}
   \right)
 \] 
is regular. Which completes the proof of the lemma.\hfill $\blacksquare$

\begin{theo}\label{fx poly caract fibré}
Let $ E \to M $ and $ F \to N $ be two vector bundles of respective ranks $ n $ and $ n '$. The $ \mathbb {R} - $ algebras $ \cl A (E) $ and $ \cl A (F), $ of polynomial functions along the fibers, are isomorphic if and only if the vector bundles $ E $ and $ F $ are isomorphic.
\end{theo}
\pre
Let $ \Psi: \cl A(F)\to \cl A(E)$ be an isomorphism of $\mathbb{R}-$algebras. Then, according to the Proposition \ref{Thomas env lisse} and the Proposition \ref{prop Jet env lisse des polynomes}, $\Psi
$ extends to an isomorphism of $\mathbb{R}-$algebras $\overline{\Psi}:{\rm C}^\infty(F)\to{\rm C}^\infty(E).$ 

Therefore, according to Milnor result, there exists a diffeomorphism $\varphi: E\to F$ such that,
 \[
   \overline{\Psi}(h)=h\circ\varphi, \, \forall h\in {\rm C}^\infty(F).
 \]
In particular, since $ \overline {\Psi} $ extends $ \Psi, $ we can write
  \begin{equation}\label{(*)}
    \cl A(E)=\{u\circ\varphi: u\in\cl A(F)\} 
  \end{equation}
In view of the Lemma \ref{lem:degree 0}, we have
 \[
   \Psi\circ\pi_F^*({\rm C}^\infty(N))=\pi_E^*({\rm C}^\infty(M)).
 \]   
We deduce the existence of a diffeomorphism $\phi:M\to N$ with
    \[  
       \Psi\circ\pi_F^*(g)=\pi_E^*(g\circ \phi),\quad \forall g\in {\rm C}^\infty(N)
    \]
Therefore, since $\pi^*_F(g)\circ \varphi=\Psi(\pi^*_F(g)),$ we have, for any $e\in E,$
\[
   g(\pi_F(\varphi(e)))=g(\phi(\pi_E(e)).
\] 
It becomes
 \begin{equation}\label{(**)}
   \pi_F\circ\varphi=\phi\circ\pi_E 
 \end{equation}
We have thus just shown that $ (\phi, \varphi) $ is an isomorphism of $ E \to M $ to $ F \to N, $ seen as differential fibrations, the local trivialization diffeomorphisms being the same as those which make vector bundles. We will therefore assume that the manifolds $ M $ and $ N $ coincide and consider $ \varphi $ as a $ M- $ isomorphism of fibrations between $E\to M$ and $F\to M.$\\

Consider now an open cover $ (U_\alpha) $ of $ M $ by trivialization domains of both $ E $ and $ F $, where $ \sigma_\alpha $ and $ \rho_\alpha $ are the diffeomorphisms of local trivialization relating to $ E $ and $ F $ respectively. We then have, for any pair of indices $ (\alpha, \beta), $ the relation
 \begin{equation}\label{(***)}
    \rho_{\beta\alpha}\circ\varphi_{\beta\alpha}=\varphi_{\beta\alpha}\circ\sigma_{\beta\alpha}
  \end{equation}  
where,  
  \[
    \varphi_{\beta\alpha}: U_{\alpha\beta}\times \mathbb{R}^n \to U_{\alpha\beta}\times\mathbb{R}^n:  (x,y)\mapsto (x,\delta_{\beta\alpha}(x,y))
  \]
is the restriction of $\varphi_{\alpha}=\rho_\alpha^{-1}\circ(\varphi|_{\pi_E^{-1}(U_\alpha)})\circ\sigma_\alpha.$ 
We further observe that  
\[  
  \varphi^x_{\alpha}: \mathbb{R}^n\to \mathbb{R}^n: y\mapsto \delta_{\alpha}(x,y)
\]
is a  diffeomorphism, polynomial in $y.$ \\

It is such that the following map
 \[
  P\in {\rm Pol }(\mathbb{R}^n)\mapsto P\circ\varphi^x_{\alpha}\in \rm{ Pol }(\mathbb{R}^n)
 \]
 is an isomorphism of $\mathbb{R}-$algebras. 
We write $ \varphi^x _{\alpha} (y) $ in the form
  \[
    \left(\lambda_{\alpha,1}^0(x)+\sum_{i=1}^{n}\lambda_{\alpha,1}^i(x)y^i+\varphi_{\alpha,1}^{\geq2}(x,y),\ldots,\lambda_{\alpha,n}^0(x)+\sum_{i=1}^{n}\lambda_{\alpha,n}^i(x)y^i+\varphi_{\alpha,n}^{\geq2}(x,y)\right)
  \] 
 where for any $i\in [1,n]\cap\mathbb{N},$  $\varphi_{\alpha,i}^{\geq2}(x,y)$ is polynomial in $y$ and each of whose terms is of degree greater than or equal to 2.\\
According to Lemma (\ref{lem:iso poly induit iso linéaire}), we observe that $(\lambda_{\alpha,j}^i(x))$ is then a regular matrix. Hence, for any index $\alpha,$ we have a diffeomorphism
\[
  \psi_\alpha:U_\alpha\times\mathbb{R}^n\to U_\alpha\times\mathbb{R}^n:(x,y)\mapsto \left(x,(\lambda_{\alpha,j}^i(x))(y)\right)
\]
 such that 
 \[
   \psi_\alpha^x:y\mapsto \left(\sum_{i=1}^{n}\lambda_{\alpha,1}^i(x)y^i,\ldots,\sum_{i=1}^{n}\lambda_{\alpha,n}^i(x)y^i\right)
 \] 
 is an automorphism of the $ \mathbb {R} -$vector space $ \mathbb {R}^n. $
Note that in virtue of the same Lemma (\ref{lem:iso poly induit iso linéaire}) we have
\[
 \psi^x_\alpha=\left(\varphi_\alpha^x\right)_{*0}\cdot
\] 
The equality (\ref{(***)}) allows therefore to write
\[
  \rho_{\beta\alpha}^x\circ\left(\psi_{\beta\alpha}^x\right)=\left(\psi_{\beta\alpha}^x\right)\circ\sigma_{\beta\alpha}^x,  
\]
by setting for any pair of indices $ (\alpha,\beta),$
 \[
   \psi_{\beta\alpha}=\psi_\alpha|_{U_{\beta\alpha}\times\mathbb{R}^n}.
 \]
Indeed, since the diffeomorphism $\rho^x: \mathbb{R}^n\to\pi_F^{-1}(x)$ is a linear one, $\pi_F^{-1}(x)$ being seen as submanifold of $F,$ we thus have 
$
   (\rho_{\beta\alpha}^x)_{*}=\rho_{\beta\alpha}^x\cdot
$
We can state an analogous result for the trivialization diffeomorphism $\sigma.$\\
We then obtain the relation 
 \[
   \rho_{\beta\alpha}\circ\psi_{\beta\alpha}=\psi_{\beta\alpha}\circ\sigma_{\beta\alpha}\cdot
 \]
We deduce from the above that there exists a unique $ M-$isomorphism of a vector bundle $ \psi: E \to F $ such that for any index $ \alpha $
 \[
 \hspace*{2cm} \psi_\alpha=\rho_\alpha^{-1}\circ\psi\circ\sigma_\alpha\cdot \qquad\qquad\qquad\qquad\qquad\blacksquare
 \] 

\begin{rem}\label{ rem induit iso gradué sur les fx poly}
  Observe  that the vector bundles isomorphism  $\psi,$ obtained through the previous Theorem \ref{fx poly caract fibré} proof, induces, by the relation $h\in{\rm C}^\infty(F)\mapsto h\circ\psi\in{\rm C}^\infty(E),$ an isomorphism of $\mathbb{R}-$algebras whose the restriction on fiberwise polynomial functions is a isomorphism of $\mathbb{R}-$algebras between $\cl A(F)$ and $\cl A(E)$ respecting the gradation.
\end{rem}

\section{Algebraic characterization of manifolds}

 The following result achieves what we start in section \ref{GraPo}.
\begin{theo}\label{caract alg variétés}
Let $M$ be a manifold.
\begin{enumerate}
\item[(a)] Any isomorphism $\Phi: \cl S(M)\to\cl S(N)$ of \,\,$\mathbb{R}-$associative algebras induces an isomorphism of\,\, $\mathbb{R}-$algebras that respects gradation.
\item[(b)] The associative \, $\mathbb{R}-$algebras  $\cl S(M)$ and $\cl S(N)$ are isomorphic if and only if $M$ and $N$ are diffeomorphic manifolds.
\end{enumerate}
\end{theo}
\pre
The statement $(a)$ is an immediate consequence of Remark \ref{ rem induit iso gradué sur les fx poly} .
 
The point $(b)$ can be obtained as corollary of the previous point $(a)$ or by direct application of  Lemme \ref{lem:degree 0}.

\section{Algebraic characterization of vector bundles}

Let $E\to M $ be a vector bundle of rank $n$ and $ T^*E\to E$ be the cotangent bundle of $E.$
We denote by $\cl D_{\cl E}(E)$ the Lie algebra generated by all differential operators $D : {\rm C}^\infty(E) \to {\rm C}^\infty(E)$ which are eigenvectors of $L_{\cl E}$, the Lie derivative in the direction of the Euler vector field.\\

Consider the associative algebra $\cl S_{\cl E}(E)=\{\sigma(T):T\in \cl D_{\cl E}(E)\},$ where $\sigma$ stands for the principal symbol operator. This space is a sub-algebra of ${\rm Pol}(T^*E):=\cl S(E),$ elements of the last space being functions that are polynomial along the fibers of $T^*E.$\\

We can now state the following result, the justification being the same as for point $(a)$ of Theorem \ref{caract alg variétés}
\begin{theo}
  Let $E\to M$ and $F\to N$ be vector bundles. 
  The associative  algebras $\cl D_{\cl E}(E)$ and $\cl D_{\cl E}(F)$ are isomorphic if and only if the vector bundles $E$ and $F$ are isomorphic\footnote{Indeed, it shown in \cite{LecLeuZih} that the Lie-algebra $\cl D_{\cl E}(E)$ characterizes the vector bundle $E$}.
\end{theo}
\begin{prop}
 The smooth envelope of the geometric $\mathbb{R}-$algebra $\cl S_{\cl E}(E)$ is given by 
 \[
    \overline{\cl S_{\cl E}(E)}={\rm C}^\infty(T^*E).
 \]
\end{prop}

\pre
We observe that
\[
  \overline{\cl S_{\cl E}(E)}\supset \pi_{T^*E}^*({\rm C}^\infty(E)).
\]
Indeed, in virtue of Proposition \ref{prop Jet env lisse des polynomes}, we have on the one hand, 
  \[
    \pi_{T^*E}^*({\rm C}^\infty(E))=\pi_{T^*E}^*(\overline{\cl A(E)})  
  \]
and on the other,
   \[
    {\rm C}^\infty(E)=\overline{\cl A(E)}\cong \overline{\pi_{T^*E}^*(\cl A(E))}
   \]  
We thus obtain the identification
   \[
     \pi_{T^*E}^*(\overline{\cl A(E)})\cong\overline{\pi_{T^*E}^*(\cl A(E))}; 
   \]   
which makes it possible to have the result announced by applying the previous Proposition \ref{inclusion egalite envpe lisse} to the inclusion    
  \[
  \pi_{T^*E}^*(\cl A(E))\subset \overline{\cl S_{\cl E}(E)}.
  \]
Let consider $\vu\in {\rm Pol}^k(T^*E).$ For any $\ve\in T^*E,$ there is a canonical chart domain $\vU\ni\ve$ of $T^*E$ associated with a trivialization $(V,\psi)$ of $E$ in which we can write 
\[
  \vu(\ve)=\sum_{|r|+|s|=k}f_{rs}(x,y)\xi^r\eta^s,
\]
where $(x,y)$ are local coordinates in $V$, $(x,y,\xi,\eta)$ being those corresponding in $\vU$ 
and $f_{rs}\in {\rm C}^\infty(V).$ 

Let us cover $ M $ by such trivialization domains and extract a Palais cover from them. 
Denote the open cover of $ T^*E $ associated with this Palais cover of $M$ by $\cl O=\cl O_1\cup\cdots\cup\cl O_N$ with $N\in\mathbb{N}$ and denote by  $\vU_{i,\alpha}(\alpha\in J)$ the elements of $\cl O_i.$  \\
Consider a partition of the unit $ \varphi_{i, \alpha} $ of $ M $ subordinate to the cover $ (U_{i,\alpha}) $ of $ M $ associated with $ \cl O.$ \\
We know that each $ \varphi_{i, \alpha} $ is compactly supported in $ U_{i, \alpha}.$ For each $ U_{i, \alpha}, $ consider a function $ \phi_ {i, \alpha} $ with compact support in $ U_ {i, \alpha} $ and which is equal to $ 1 $ where $\varphi_{i,\alpha}$  is nonvanishing.

We then have 
 \begin{eqnarray*}
\pi_{T^*E}^*\circ\pi_E^*(\varphi_{i,\alpha})\vu & = & \pi_{T^*E}^*\circ\pi_E^*(\varphi_{i,   
                                                                   \alpha}\cdot\phi_{i,\alpha})\vu\\
                                                & = & \sum_{|r|+|s|=k} g^{rs}_{i,\alpha}\cdot\vv^{rs}_{i,\alpha}.
 \end{eqnarray*}
In the above equalities, we have that $g^{rs}_{i,\alpha}\in \pi_{T^*E}^*({\rm C}^\infty(E))$ and $\vv^{rs}_{i,\alpha} \in \cl S_{\cl E}(E)$ are of support in $\vU_{i,\alpha}$ and are given by
  \[
    g^{rs}_{i,\alpha}=\varphi_{i,\alpha}f^{rs}_{i,\alpha} \mbox{ and } \vv^{rs}_{i,\alpha}=\phi_{i,\alpha}\xi^r\eta^s.
  \]
Therefore, we can write
\begin{eqnarray*}
  \vu & = &  \sum_{|r|+|s|=k}\sum_{i,\alpha}g^{rs}_{i,\alpha}\cdot \vv^{rs}_{i,\alpha}\\
   & = & \sum_{|r|+|s|=k}\sum_{i}\left(\sum_{\alpha}g^{rs}_{i,\alpha}\right)\cdot \left(\sum_{\beta}\vv^{rs}_{i,\beta}\right)\\
    & = &\sum_{|r|+|s|=k}\sum_{i} g^{rs}_{i}\cdot \vv^{rs}_{i}
\end{eqnarray*}
since for $\alpha\neq\beta,$ one has $g^{rs}_{i_\alpha}\cdot\vv^{rs}_{i_\beta}=0.$\\
We then obtain $\vu\in\overline{\cl S_{\cl E}(E)},$ since we have just shown that $ \vu $ decomposes into a sum whose terms are products having two factors, one in $\pi_{T^*E}^*({\rm C}^\infty(E))$ and the other in $\cl S_{\cl E}(E).$ Therefore,
  \[
   {\rm Pol}(T^*E)\subset \overline{\cl S_{\cl E}(E)}.
  \]
According to previous Proposition \ref{inclusion egalite envpe lisse}, we thus have the inclusion
 \[
   \overline{ {\rm Pol}(T^*E)}\subset\overline{\cl S_{\cl E}(E)}.
 \]  
But we also have $\overline{{\rm Pol}(T^*E)}={\rm C}^\infty(T^*E)$. It then comes  
\[
 {\rm C}^\infty(T^*E) \subset\overline{\cl S_{\cl E}(E)}.
\]
 We then deduce the equality 
\[
 \overline{\cl S_{\cl E}(E)}={\rm C}^\infty(T^*E);
\]
because the inclusion $\overline{\cl S_{\cl E}(E)}\subset {\rm C}^\infty(T^*E)$ comes from the Proposition \ref{inclusion egalite envpe lisse}, applied to the following inclusions  
\[
\hspace*{2cm}\cl S_{\cl E}(E)\subset {\rm Pol}(T^*E)\subset\overline{{\rm Pol}(T^*E)}.\hspace*{3cm} \blacksquare
\]

We can now state the following result which will allow us to draw some conclusions on vector bundles, in relation to the structure of $ \mathbb {R} - $ algebra of the space of symbols of homogeneous operators.

\begin{prop}\label{iso sde vers sdf  indui iso poly}
Let $E\to M$ and $F \to N$ be two vector bundles and $T^*E\to E$ and $T^*F\to F$ be the cotangent bundles associated with $E$ and $F$ respectively.
Any isomorphism of $\mathbb{R}-$algebras $\Psi:\cl S_{\cl E}(E)\to \cl S_{\cl E}(F)$ extends to an unique isomorphism of $\mathbb{R}-$algebras $\overline{\Psi}:{\rm C}^\infty(T^*E)\to {\rm C}^\infty(T^*F)$ such that
 \[
  \overline{\Psi}(\cl A^0((E)))=\cl A^0((F))  
 \]
where $\cl A^0((E))=\pi_{T^*E}^{*}(\cl A^0(E))$ and $\cl A^0((F))$ is the $\mathbb{R}-$algebra defined analogously. 
\end{prop}
\pre
If $\vu\in \cl A^0((E))\subset {\rm Pol}^0(T^*E),$ is nonvanishing on $T^*E,$ it is then the same for 
  \[
  \vu^{-1}:T^*E\to \mathbb{R}:\ve\mapsto \frac{1}{\vu(\ve)}
  \]
 which  is an element of ${\rm Pol}^0(T^*E).$ But since $\vu$ is a  polynomial function of zero degree in $y,$ when we consider a system of local coordinates $(x,y,\xi,\eta),$ in a canonical chart of $T^*E,$ the same is true for $\vu^{-1};$ and thus  
 \[
 \vu^{-1}\in\cl A^0((E))\subset\cl S^0_{\cl E}(E)= \pi_{T^*E}^{*}(\cl A(E)).
 \]
Therefore, it is obvious that we have 
 \[
\overline{\Psi}(\vu)\cdot\overline{\Psi}(\vu^{-1})=1_{T^*F}:\ve\mapsto 1 \cdot 
 \]
This becomes, because $\overline{\Psi}$ extends $\Psi,$
\[
\Psi(\vu)\cdot\Psi(\vu^{-1})=1_{T^*F}
\]
The above equality allows to conclude that  $\Psi(\vu)$ and $\Psi(\vu^{-1})$ are constant along the fibers of $T^*F.$ Thus, $\Psi(\vu)$ and $\Psi(\vu^{-1})$ are elements of 
\[
\cl S^0_{\cl E}(F)=\pi_{T^*F}^*(\cl A(F)).
\]  
But since invertible elements of $\cl A(F)$ that have their inverses in $\cl A(F)$ are in $\cl A^0(F)$ and 
\[
\pi_{T^*F}^{*}|_{{\rm C}^\infty(F)}:{\rm C}^\infty(F)\to {\rm C}^\infty(F)
\]
is an isomorphism of $\mathbb{R}-$algebras, we also conclude that $\Psi(\vu)$ is in $\cl A^0((F)).$ 
Since for any $\vu\in\cl A^0((E)),$ the element $\vu^2+1$ is a nonvanishing one, we have that  $(\Psi(\vu))^2$ is an element of $\cl A^0((F)).$ 
We directly deduce that 
\[
\Psi(\vu)\in\cl A^0((F)).
\]
And the above relation achieves the proof of the proposition. \hfill $\blacksquare$

\begin{rem}
 Let $E\to M$ be a vector bundle of rank $n$ with ${\rm dim }(M)=m.$\\ To any transition diffeomorphism of the vector bundle $E\to M$ of the form 
  $(x,y)\mapsto (x,A(x)(y))$ corresponds a transition diffeomorphism of the tangent bundle $TE\to E$ of the form
  \[
     (x,y,\vh,\vk)\mapsto (x,y,\vh,(A'(x)\cdot y)(\vh)+A(x)(\vk)),  
  \]
 where $A'(x)\cdot y$ is a $(n,m)-$type matrix  such that
  \[
    (A'(x)\cdot y)(\vh)=(A_{*_x}\vh)(y).
  \]
\end{rem}
The transition diffeomorphisms of $T^*E\to E$ are then given by
 \[
   (x,y,\xi,\eta)\mapsto (x,y, \xi-^t(A^{-1}(x)\circ(A'(x)\cdot y))(\eta)\, ,\, ^tA^{-1}(x)(\eta))
 \]
These diffeomorphisms thus define  a differential fibration $T^*E\to M$ whose  projection is given by 
  \[
   \pi: T^*E\to M: \ve\mapsto \pi_{E}\circ\pi_{T^*E}(\ve)
  \]
\begin{cor}
Let $E\to M$ and $F\to N$ be vector bundles. If the $\mathbb{R}-$algebras $\cl S_{\cl E}(E)$ and $\cl S_{\cl E}(F)$ are isomorphic, then the differential fibrations $T^*E\to M$ and $T^*F\to N$ are isomorphic. 
\end{cor}
\pre
Indeed, in virtue of Proposition \ref{iso sde vers sdf  indui iso poly}, the $\mathbb{R}-$algebras ${\rm C}^\infty(T^*E)$ and ${\rm C}^\infty(T^*F)$ are isomorphic.\\

There then exists a diffeomorphism
 $\Phi:T^*E\to T^*F$ such that, in accordance with the proposition cited above, the $\mathbb{R}-$algebras isomorphism  $\Psi$ between $\cl S_{\cl E}(F)$ and $\cl S_{\cl E}(E)$  is given by
 \[
   \Psi: \vu\in \cl S_{\cl E}(F) \mapsto \vu\circ\Phi\in \cl S_{\cl E}(E)\cdot
 \]
The same Proposition \ref{iso sde vers sdf  indui iso poly} allows to write
 \[
   \Psi(\pi_{T^*F}^*\circ\pi_F^*({\rm C}^\infty(N)))=\pi_{T^*E}^*\circ\pi_E^*({\rm C}^\infty(M))\cdot
 \] 
Therefore, the restriction of $\Psi$ to $\cl A^0(F)$ induces an isomorphism of $\mathbb{R}$-algebras $\underline{\Psi}$ between ${\rm C}^\infty(N)$ and ${\rm C}^\infty(M).$
 
By Milnor result, there exists a diffeomorphism $\phi: M\to N$ such that 
 \[
 \underline{\Psi}(g)=g\circ\phi.
 \]
We deduce from the above,
 \[
   \Psi(\pi_{T^*F}^*\circ\pi_F^*(g))=\pi_{T^*E}^*\circ\pi_E^*(g\circ\phi)\cdot
 \]
But we have the equality
 \[
   (\pi_{T^*F}^*\circ\pi_F^*(g))\circ\Phi=\Psi(\pi_{T^*F}^*\circ\pi_F^*(g))\cdot
 \]
Hence, for any  $\ve\in T^*E,$ we have
 \[
  (\pi_{T^*F}^*\circ\pi_F^*(g))(\Phi(\ve))=\pi_{T^*E}^*\circ\pi_E^*(g\circ\phi)(\ve)\cdot
 \]
This is equivalent to
 \[
   g(\pi_F\circ\pi_{T^*F}(\Phi(\ve)))=g(\phi(\pi_E\circ\pi_{T^*E}(\ve))), \forall \ve\in T^*E, \forall g\in{\rm C}^\infty(N)\cdot
 \]
We deduce that
   \[
     (\pi_F\circ\pi_{T^*F})\circ\Phi=\phi\circ(\pi_E\circ\pi_{T^*E});
   \]
and the announced result is proved. \hfill $\blacksquare$\\


The above corollary also allows us to say that if the $  \mathbb {R} - $ algebras $ \cl S_{\cl E} (E) $ and $ \cl S_{\cl E} (F) $ are isomorphic, the basis $ M $ and $ N $ of the vector bundles $ E $ and $ F $ are diffeomorphic. It also implies that these bundles have the same rank. \\

To obtain a characterization of vector bundles, let us return to the previous $ \Phi $ isomorphism and keep the notations of the proof of Theorem \ref{fx poly caract fibré}.\\

We can therefore consider $ \Phi $ as an $ M-$isomorphism of the fibrations $ T^*E \to M $ and $ T^*F \to M $ \\
Let now $ (U_\alpha) $ be an open cover  of $ M $ by trivialization domains of $ T^*E $ and $ T^*F $. Let us denote by $ \sigma_\alpha $ and $ \rho_\alpha $ the local trivialization diffeomorphisms relating to these vector bundles.\\
We now can write
\[
   \Phi_\alpha^x: \mathbb{R}^{m+2n} \to \mathbb{R}^{m+2n}: (y,\xi, \eta)\mapsto \Delta_\alpha(x,y,\xi, \eta)
\] 
with $\Phi_{\alpha}=\rho_\alpha^{-1}\circ(\Phi|_{\pi^{-1}(U_\alpha)})\circ\sigma_\alpha.$ 

In addition, observe that $\Phi_\alpha^x$  is a diffeomorphism, polynomial on $y,\xi,\eta ,$ such that
 \[
  P\in {\rm Pol}(\mathbb{R}^{m+2n})\mapsto P\circ\Phi^x_{\alpha}\in {\rm Pol}(\mathbb{R}^{m+2n})
 \]
 is an automorphism of $\mathbb{R}-$algebra.\\
By Lemma \ref{lem:iso poly induit iso linéaire}, the linear part of $ \Phi_\alpha^x $ allows to define, for any index $ \alpha, $ a diffeomorphism 
\[
  \Psi_\alpha:U_\alpha\times\mathbb{R}^{m+2n}\to U_\alpha\times\mathbb{R}^{m+2n}
\]
We deduce, as before, that there exists a $ M- $ isomorphism of vector bundles and only one $\Psi:T^*E\to T^*F$ such that for any index $\alpha$
 \[
  \Psi_\alpha=\rho_\alpha^{-1}\circ\Psi\circ\sigma_\alpha.
 \] 
We have thus just proved the following result.
\begin{prop}
Let $E\to M$ be a vector bundle.\\ 
Then, the associative algebra $ S_{\cl E} (E) $ characterizes the vector bundle $T^*E\to M.$
\end{prop}

If rather, we define 
\[
   \Phi_\alpha^{x,y}: \mathbb{R}^{m+n} \to \mathbb{R}^{m+n}: (\xi, \eta)\mapsto \Delta_\alpha(x,y,\xi, \eta)
\] 
with $\Phi_{\alpha}=\rho_\alpha^{-1}\circ(\Phi|_{\pi^{-1}(U_\alpha)})\circ\sigma_\alpha,$ where this time, $\rho$ and $\sigma$ are seen as transition diffeomorphisms of $E\to M$ and $F\to M,$ we obtain, in a similar way as in the previous lines, a $ M- $ isomorphism of vector bundles $\Psi$ between $E\to M$ and $F\to M.$\\
We summarize this in the following statement.
\begin{theorem}
Let $E\to M$ be a vector bundle.\\ Then, the associative algebra $ S_{\cl E} (E) $ characterizes the vector bundle $E\to M.$
\end{theorem}
\section{Classical limit of homogeneous operators of zero weight}
Let $E\to M$ be a vector bundle. We denote by
 \[
 \cl D_0(E)=\{T\in\cl D(E):[\cl E,D]=0\}
 \]
the quantum Poisson algebra of homogeneous operators of zero weight.\\

Consider the classical Poisson algebra $ \cl S_0 (E) $, classical limit of $ \cl D_0 (E). $ We then have
\[
  \cl S_0(E)=\bigoplus_{k\geq0} \cl S_0^k(E),
\]
with $\cl S_0^k(E)=\{\sigma(T):T\in\cl D_0^k(E)\}.$ Seen as associative subalgebra of ${\rm Pol}(T^*E),$ we have
 \[
   \cl S_0^0(E)=\pi_{T^*E}^*(\cl A^0(E))=\cl A^0((E)).
 \] 
Locally, in a canonical chart of $ T^*E $ associated with an adapted chart of $ E, $ an element $ \vu $ of $ \cl S_0^k (E) $ is written in the form
 \begin{equation}\label{(*v)}
    \vu(\ve)=\sum_{|r|+|s|\leq k}f_{rs}(x)y^r\xi^s\eta^r, 
 \end{equation}   
where $\ve\in T^*E$ admits $(x,y,\xi,\eta)$ as local coordinates \footnote{This comes from the local form of homogeneous operators established in\cite{LecLeuZih}}.

\begin{prop}\label{iso szeroe induit iso azeroe}
   Let $E\to M$ and $F\to N$ be two vector bundles. \\
   If $\Psi:\cl S_0(E)\to \cl S_0(F)$ is an isomorphism of $\mathbb{R}-$algebras, then we have 
     \[
       \Psi(\cl A^0((E)))=\cl A^0((F)).
     \] 
\end{prop}
\pre 
The proof is analogous to that of the Lemma \ref{lem:degree 0}. \\Indeed, let $\vu\in\cl A^0((E))$ be nonvanishing on $T^*E.$ Then $\vu=\pi_{T^*E}^*(u),$ with $u\in\cl A^0(E)$ a nonvanishing function on $E.$ 
Therefore, $u^{-1}\in\cl A^0(E)$ and, using $\vu^{-1}=\pi_{T^*E}^*(u^{-1}),$ we obtain the relation 
  \[
  \vu\cdot\vu^{-1}=1_{T^*E}:\ve\mapsto 1.
  \]
Hence, we have
  \[
    \Psi(\vu)\cdot\Psi(\vu^{-1})=1_{T^*F}.
  \]
We deduce that   
  \[
  \Psi(\vu)\in {\rm Pol}^0(T^*F)\cap\cl S_0(F)=\cl A^0((F)).
  \]
Considering any element $\vu$ of $\cl A^0((E)),$ we have that $\vu^2+1$ is nonvanishing on $T^*E$, and it comes 
  $\Psi(\vu)\cdot\Psi(\vu)\in\cl A^0((F)).$ \\  
We can then conclude because it implies 
  \[
 \hspace*{4cm} \Psi(\vu)\in\cl A^0((F)) \hspace*{6cm}\blacksquare
  \]

We then deduce the following result.
\begin{cor}
 Let $ E \to M $ and $ F \to N $ be two vector bundles. If the $ \mathbb {R}-$algebras $ \cl S_0 (E) $ and $ \cl S_0 (F) $ are isomorphic then the differential manifolds $ M $ and $ N $ are diffeomorphic.
\end{cor}




\newpage

\nocite{*}
\begin{thebibliography}{ }
\addcontentsline{toc}{chapter}{Bibliographie}


\bibitem{GelKol}
Gel'fand I, Kolmogoroff A, \textit{ On rings of continuous functions on topological spaces},  C. R. (Dokl.) Acad. Sci. URSS, 22 (1939), pp. 11-15 
\bibitem{GraPon2}
 Grabowski J, Poncin N, \textit{ Automorphisms of quantum and classical Poisson algebras}, Comp. Math., \textbf{140} (2004), pp. 511-527
\bibitem{GraPon1}
  Grabowski J, Poncin N, \textit{ Lie-algebraic characterizations of manifolds}, Central Europ. J. of Math., \textbf{2(5)} (2005), pp. 811-825 
\bibitem{GraPon4}
 Grabowski J, Poncin N, \textit{ On quantum and classical Poisson algebras}, Banach Center Publ. \textbf{76}, Warszawa (2007), pp. 313-324.
\bibitem{Lec1}
Lecomte P, \textit{ On the infinitesimal automorphisms of a vector bundle}, J. Math. pure et appl. Go (1981), pp. 229-239
\bibitem{Lec2} 
Lecomte P, \textit{ On some sequence of graded Lie algebras associated to
manifolds}, Ann. Glob. Anal. Geom., 12 (1994), pp. 183-192
\bibitem{Lec3}
 Lecomte P.B.A, \textit{ Note on the Linear Endomorphisms of a Vector Bundle.}, Manuscripta mathematica 32 (1980): 231-238
\bibitem{Lec4} 
Lecomte P, \textit{ Sur l'algèbre de Lie des sections d'un fibré en algèbres de Lie, } Ann. Inst. Fourier, XXX, Fasc. \textbf{4}, (1980), pp. 35-50.
\bibitem{LecLeuZih}
 Lecomte P.B.A, Leuther T, Zihindula Mushengezi E, \textit{ On a Lie algebraic characterization of vector bundles, }SIGMA 8 (2012), 004:10 pages, 2012.
\bibitem{LecZih1}
 Lecomte P.B.A, Zihindula Mushengezi E, \textit{ On quasi quantum Poisson algebras : Lie-algebraic characterization, } arXiv:2007.14649v1 [math.DG].
\bibitem{LecZih2}
 Lecomte P.B.A, Zihindula Mushengezi E, \textit{ Lie algebra of homogeneous operators of a vector bundle, } arXiv:2007.14692v1 [math.DG]. 
 \bibitem{Leu}
Leuther T, \textit{Affine bundles are affine spaces over modules} arXiv:1201.5812v1 [math.DG]
\bibitem{Jet}
 Nestruev J, \textit{ Smooth manifolds and observables} volume 220 of Graduate Texts in Mathematics. Springer-Verlag, New York, 2003. Joint work of A. M. Astashov, A. B. Bocharov,
S. V. Duzhin, A. B. Sossinsky, A. M. Vinogradov and M. M. Vinogradov, Translated from
the 2000 Russian edition by Sossinsky, I. S. Krasil'schik and Duzhin
\bibitem{PurShan}
 Pursell L E, Shanks M E, \textit{ The Lie algebra of a smooth manifold}, Proc. Amer. Math. Soc. \textbf{5} (1954), pp. 468.
 
 


\end{thebibliography}
\end{document}